\documentclass[11pt, a4paper]{article}

\usepackage[utf8]{inputenc}
\usepackage{alphalph, amsfonts, amsmath, amssymb, booktabs, etoolbox, fancyhdr, hyperref, lipsum, mathtools, pgf, pgfplots, psfrag, scalerel, siunitx, subcaption, tikz, xcolor,url}
 \usetikzlibrary{calc}

\usepackage[margin=0.85in]{geometry}

\newtheorem{thrm}{Theorem}[section]

\newtheorem{proposition}[thrm]{Proposition}

\newtheorem{qstn}[thrm]{Question}
\newtheorem{rmrk}[thrm]{Remark}

\newenvironment{proof}{\noindent {\bf Proof. }}{\hfill $\Box$ \newline\par}

\newcommand{\Pmax}{P_{\max}}
\newcommand{\vmax}{v^{\max}}
\newcommand{\wmax}{w^{\max}}
\newcommand{\winit}{w_{\text{init}}}

\usepackage{bbm}

\title{Speed planning by minimizing travel time and energy consumption}
\author{Stefano Ardizzoni, Luca Consolini, Mattia Laurini, Marco Locatelli  % <-this % stops a space
\thanks{All authors are with the University of Parma, Department of Engineering and Architecture, Parco Area delle Scienze 181/A, 43124 Parma, Italy. E-mails:
{\tt\footnotesize \{stefano.ardizzoni, luca.consolini, mattia.laurini, marco.locatelli\}@unipr.it}
}}
\date{}
\begin{document}
\maketitle
\begin{abstract}
In this paper we address the speed planning problem for a vehicle over an assigned path with the aim of minimizing a weighted sum of travel time and energy consumption under suitable constraints (maximum allowed speed, maximum traction or braking force, maximum power consumption). The resulting mathematical model is a non--convex optimization problem. We prove that, under some mild assumptions, a convex reformulation of the non--convex problem is exact. In particular, the convex reformulation is a Second Order Cone Programming (SOCP) problem, for which efficient solvers exist. Through the numerical experiments we confirm that the convex relaxation can be solved very efficiently and, moreover, we also provide the Pareto front of the trade-off between the two objectives (travel time and energy consumption).
\end{abstract}
\section{Introduction}

The optimization of energy consumption is of key importance in road vehicles management, and in transportation in general. In this scenario, one relevant problem is the optimization of the speed of a road vehicle on an assigned path. In this paper, we present an efficient method for minimizing a weighted average of travel time and energy consumption.
% In this work, we address the problem of minimizing both travel time and power consumption for a vehicle moving along an assigned path.

\tikzset{
  auto/.pic = {
    \begin{scope}[shift={(-4,+0.1)}]
    \shade[top color=red, bottom color=white, shading angle={135}]
        [draw=black,fill=red!20,rounded corners=1.2ex,very thick] (1.5,.5) -- ++(0,1.3) --  ++(4.5,0) -- ++(1.5,-0.5) -- ++(0,-0.8) -- (1.5,.5) -- cycle;
    \draw[very thick, rounded corners=0.5ex,fill=black!20!blue!20!white,thick]  (2.5,1.8) -- ++(0.5,0.7) -- ++(1.8,0) -- ++(1,-0.7) -- (2.5,1.8);
    \draw[thick]  (3.5,1.8) -- (3.5,2.5);
    \draw[draw=black,fill=gray!50,thick] (2.75,.5) circle (.6);
    \draw[draw=black,fill=gray!50,thick] (5.5,.5) circle (.6);
    \draw[draw=black,fill=gray!80,semithick] (2.75,.5) circle (.5);
    \draw[draw=black,fill=gray!80,semithick] (5.5,.5) circle (.5);
    \end{scope}
  }
}

\begin{figure}
  \begin{center}
\begin{tikzpicture}[scale=1.5,domain=-1:6]
      \draw[color=green,thick] plot (\x,{0.1*(\x*\x-5*\x)})  node[right] {};%{$\gamma(s)$};
\draw (4,{0.1*(4*4-5*4)})  node(c) {};  
\pic at (c)  [draw,scale=0.5,rotate=18] {auto};
\draw[->] (4,{0.1*(4*4-5*4)})  -- ++ (1.8,{0.18*(8-5});  
\draw[dashed] (4,{0.1*(4*4-5*4)})  -- ++ (1.5,0) node (a){} -- ++ (1.0,0);
\draw [-latex] (a) arc (0:{atan(0.3)}:1.5) node [right,midway] {$\alpha(s)$};
\draw  [fill] (4,{0.1*(4*4-5*4)}) node {}  circle(0.05);
\draw  (4,{0.1*(4*4-5*4)}) node [below]  {$\gamma(s)$};

\end{tikzpicture}
\end{center}
\caption{Problem setting.}
\label{fig_initial}
\end{figure}
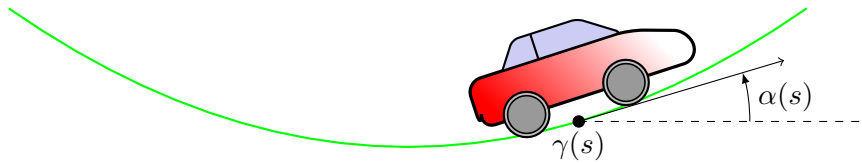

Figure~\ref{fig_initial} shows the basic problem setting. We consider a vehicle on a road with a known elevation profile $\gamma(s)$ and slope $\alpha(s)$, where $s$ is the arc-length position on the travelled curve. We want to find the speed law that optimizes a weighted sum of the travel time and the overall used energy. We take into account maximum speed constraints along the path, forces dissipated by aerodynamic effects and tyre rolling resistance, no slipping constraints, and constraints on the maximum available power.
In more detail, we consider the following optimization problem, which is the discretized version of a continuous-time one, that is derived in detail in Appendix~\ref{sec:prob_deriv}:
\begin{align}
& \min_{w,F} \sum_{i=1}^{n-1} h \left(\lambda \max\{\eta F_i,F_i\} + \frac{1}{\sqrt{w_i}}\right) & \label{obj_fcn} \\
\textrm{such that} & & \nonumber \\
& w_i \leq \wmax_i & \text{for } i \in \{1, \ldots, n\}, \label{constr:max_speed} \\
& |F_i| \leq Mg\mu & \text{for } i \in \{1, \ldots, n-1\}, \label{constr:max_force} \\
& F_i\sqrt{w_i} \leq \Pmax & \text{for } i \in \{1, \ldots, n-1\}, \label{constr:max_power} \\
& \frac{M}{h} (w_{i+1}-w_i) = - \Gamma w_i + F_i - M g \sin \alpha_i - M g c & \text{for } i \in \{1, \ldots, n-1\}, \label{constr:forces} \\
& w_i \geq 0 & \text{for } i \in \{1, \ldots, n\}, \label{constr:min_speed} \\
& w_1 = \winit. & \label{constr:init_cond}
\end{align}
Here, $w_i$ and $F_i$ are the optimization variables. They represent the squared speed and the traction (or braking) force at position $i$ along a given path, for $i \in \{1, \ldots, n-1\}$.
Constraints~\eqref{constr:max_speed} require that the squared speed be no larger than a given maximum squared speed profile $\wmax_i$, for $i \in \{1, \ldots, n\}$. Note that this constraint also allows bounding maximum lateral accelaration (see our previous work~\cite{consolini2017scl}).
Constraints~\eqref{constr:max_force} represent the maximum traction or braking force, where $\mu$ is the road/tyre friction coefficient and $g$ the gravitational acceleration.
Constraints~\eqref{constr:max_power} model maximum power consumption, with $\Pmax$ being the maximum available power. These constraints come into play only when force $F_i$ is positive, that is when the vehicle is accelerating, and there is a traction force acting on it (i.e., $F_i > 0$) and not a braking one (i.e., $F_i \leq 0$).
Constraint~\eqref{constr:forces} represents the vehicle longitudinal dynamic equation. By Newton's law,  mass ($M$) times acceleration (approximated through a finite difference with length step $h$) is equal to the sum of the forces acting on the vehicle. These forces are the aerodynamic drag force $-\Gamma w_i$, the traction (or braking) force $F_i$, the force $-Mg\sin(\alpha_i)$ due to road grade $\alpha_i$, that depends on position $i$ along the path, and the tires rolling resistance $- M g c$.
Finally, constraints~\eqref{constr:min_speed} require that the squared speed be nonnegative, whilst~\eqref{constr:init_cond} imposes the initial squared speed $\winit \geq 0$.
The objective function~\eqref{obj_fcn} represents the weighted sum of the consumed energy and the overall travel time. Constant $\eta \in [0, 1]$ represents the relative amount of energy that can be recovered by regenerative braking. Parameter $\eta$ can be significantly large (for instance $0.7$) in electric engines. Smaller positive values $\eta$ can represent a hybrid engine, since the maximum regenerative capacity of braking is smaller. In thermal engines we can set $\eta = 0$, since all kinetic energy is dissipated by braking.
Note that the objective term $\sum_{i=1}^{n-1} h \left(\lambda \max\{\eta F_i,F_i\}\right)$ only represents the overall mechanical energy in traction. We do not take into account energy dissipated in other forms. In particular, in thermal engines, a large fraction of the overall energy is dissipated in heat, which we do not account for.
Note that $\lambda$ allows weighting power consumption and travel time in the objective function. Constant $\lambda$ represents the time cost (in seconds) per Joule (its dimension is \unit{\second\per\joule}).
Problem~(\ref{obj_fcn})--(\ref{constr:init_cond}) is non--convex due to power constraint~\eqref{constr:max_power}.
%The problem is modeled as a constrained optimization one, involving the speed and force of the vehicle, subject to physical and operational limits such as maximum speed, power constraints, and dynamic forces.
Nonconvexity makes this problem challenging to solve efficiently and exactly, motivating the use of reformulations and relaxations to both improve computational performances, and obtain global optimality results on the solutions.
%The problem of speed planning can be approached following two main paradigms:
%\begin{enumerate}
%\item Simultaneous trajectory optimization, where the path and its timing law a%re jointly designed.
%\item Sequential path and speed optimization, where the geometric path is first% determined, followed by the optimization of the speed law along such path.
%\end{enumerate}
%The former approach guarantees optimality but is computationally demanding.
%The second approach is sub-optimal, but allows for separate optimization steps,% making it more practical for specific applications.
%In this work, we focus on the second approach, assuming that we are given a predefined path.
% Of course, these two objectives are competing with each other and the way they are optimized influences a vehicle performances both in terms of travel time and energy consumption.
%This can be particularly relevant for electric vehicles and autonomous driving, for which energy efficiency can directly impact effectiveness.
\subsection{Literature Review} 

Problem~(\ref{obj_fcn})--(\ref{constr:init_cond}) becomes a minimum-time one if $\lambda=0$. On the other hand, if $\lambda$ is large, this problem becomes a minimum-fuel one. We review the literature both on minimum-time and minimum-fuel problems, since they are strongly related.
\subsubsection{Minimum-time speed planning}
Various works focus on the minimum-time control problem for road vehicles on an assigned path.
%These works considers different models, or different kinematic, dynamic, and mechanical constraints. For instance, these constraints can be velocities, accelerations and maximal steering
%angle (see, for instance,~\cite{bianco2006optimal,Frego16,1700047,doi:10.1177/027836498600500304,velenisMinimumTimeTravelVehicle2008}).
Some works represent the vehicle speed as a function of time, others represent the vehicle speed as a function of the arc-length position along the assigned path. In some cases, the second choice considerably simplifies the resulting optimization problem.
Among the works that make the first choice, \cite{MunOllPraSim:94,munoz1998speed} propose an iterative algorithm based on the concatenation of third degree polynomials. Reference~\cite{SolNun2006} achieves a solution with an algorithm based on the five-splines scheme of~\cite{bianco2006optimal}, while~\cite{Villagra-et-al2012} proposes a closed-form speed profile obtained as a concatenation of three different classes of speed profiles. Finally,~\cite{CheHeBuHanZha2014} proposes an algorithm which returns a piecewise linear speed profile.
%Other works (see, for instance, represent the speed law as a function $v$ of the arc-length position $s$ and not as a function of time.
Among the works that parameterize the speed as a function of the arc-lenght,~\cite{Minari16,minSCL17} consider a problem that, after finite-element discretization, can be reformulated as a convex one. This allows deriving efficient solution methods such as those presented in~\cite{CLMNV19,minSCL17,LippBoyd2014}. In particular, reference~\cite{minSCL17} presents an algorithm, with linear-time computational complexity with respect to the number of variables, that provides an optimal solution after spatial discretization. Namely, the path is divided into $n$ intervals of equal length, and the problem is approximated with a finite dimensional one in which the derivative of $v$ is substituted with a finite difference approximation. Some other works directly solve this problem in continuous-time~\cite{consolini2020solution,Frego16,Velenis2008}. In particular, paper~\cite{consolini2020solution} computes directly the exact continuous-time solution without performing a finite-dimensional reduction, proposing a method which is very simple and can be implemented very efficiently, so that it is particularly well suited for real-time speed planning applications. Reference~\cite{RaiPerCGL17jerk} solves the problem with an additional jerk constraint through a heuristic approach that computes a speed profile by bisection; the proposed method is very efficient, however, the optimality of the obtained profile is not guaranteed.

%As an alternative to minimize travel time, some papers (see, e.g.,~\cite{baum_speed-consumption_2014,kim_real-time_2020,ozatay_cloud-based_2014}) also take into account energy consumption. Some studies focus on finding the path that allows the vehicle to consume the least amount of energy.
%In particular, various work study the problem of computing routes for electric vehicles in road networks~\cite{baum2013energy,eisner2011optimal,sachenbacher2011efficient,storandt2012quick,baum2014speed}. Since their
%battery capacity is limited, and consumed energy per distance increases with velocity, driving
%the fastest route is often not desirable, and may even be infeasible. On the other hand, the
%energy-optimal route may be too conservative in that it contains unnecessary detours or simply
%takes too long.
%Most papers have focused on the integration of battery capacity
%constraints into classical single-criterion
%route planning algorithms optimizing energy consumption~\cite{baum2013energy,eisner2011optimal,sachenbacher2011efficient}. However, such routes
%may have disproportionate detours: driving slower saves energy at the cost of greatly longer
%travel time. Storandt~\cite{storandt2012quick} therefore optimizes energy consumption, but bounding the amount
%by which travel time may increase.

%Reference~\cite{baum2014speed}  proposes to use multicriteria optimization to obtain Pareto sets
%of routes that trade energy consumption for speed. 

\subsubsection{Fuel consumption minimization}

Other works seek the speed profile that minimizes fuel consumption on an assigned path.
Various works optimize the energy utilization of metro vehicles for a given trip. For instance,~\cite{wang2011optimal} formulates the optimal speed profile for
 energy saving problem as a Mixed Integer
 Linear Programming (MILP) one. In~\cite{kang2011ga}, Kang et
 al. present an algorithm for optimizing a train speed
 profile by controlling the
 coasting point, using Genetic Algorithms (GA).
 Reference~\cite{calderaro2014algorithm} finds the sequence of basic control regimes that
 minimizes the energy consumed on a given path, taking into
 account the track topology (slopes and curves), the
 mechanical characteristics of the vehicle, the electrical
 characteristics of the feeding line, and the effect of
 regenerative braking. 
\newline\newline\noindent
 Other works minimize fuel consumption for private vehicles.
However, human drivers are not aware of the optimal velocity profile for a given
 route. Indeed, the globally optimal velocity trajectory depends on
 many factors, and its calculation requires intensive computations~\cite{gustafsson2009automotive,russell2002integrated}. With advancements in communication, sensors, and in-vehicle computing, real-time optimizations are becoming more feasible. For example, some public transportation vehicles in Europe now communicate with traffic lights~\cite{koenders2008cooperative}, and the USA is experimenting with broadcasting red light timings for safety~\cite{intersections2008cooperative}. Several algorithms have been proposed for optimizing speed profiles, such as predicting optimal velocity profiles for approaching traffic lights~\cite{asadi2010predictive}, or minimizing energy use for a given route with a traffic light~\cite{ozatay2012analytical}. Despite these developments, many algorithms still require expensive on-board systems and have limited real-time applicability. To address these limitations, cloud computing offers a promising solution for real-time speed optimization~\cite{wollaeger2012cloud}. Reference~\cite{ozatay2014cloud} builds on that by implementing a real-time Speed Advisory System (SAS), that uses cloud computing to generate optimal velocity profiles based on traffic and geographic data, offering a global optimization approach. 
 A more aggressive approach is rapidly accelerating to a
 given speed, followed by a period of coasting to a lower speed,
 which is called the “pulse-and-glide” (PnG) strategy~\cite{li2015effect,li2016fuel,li2015mechanism,xu2015fuel}. PnG driving can
 achieve significant fuel savings in vehicles with continuously
 variable transmissions (CVT)~\cite{li2015mechanism} and step-gear mechanical
 transmissions~\cite{xu2015fuel}. The optimality and the increasing efficiency are proven by theoretical calculations. Reference~\cite{kim2019real} introduces a real-time implementable pulse-and-glide algorithm in the speed-acceleration domain, achieving 3\%–5\% fuel savings without compromising optimality. Some papers focus on the study of road slope~\cite{hellstrom2009look,kamal2011ecological,musardo2005ecms,xu2018design}, which significantly impacts fuel economy, requiring optimization of speed and control systems. Model predictive control (MPC) methods, like those by Kamal et al.~\cite{kamal2011ecological} and Hellström et al.~\cite{hellstrom2009look}, improve fuel efficiency, but face high computational demands due to system nonlinearity. Non-predictive methods, such as ECMS~\cite{musardo2005ecms}, offer reduced computational loads and near-optimal fuel savings by relying on instantaneous slope data. Both approaches are viable for connected automated vehicles (CAVs), but require balancing fuel economy and computational efficiency for practical implementation.
\newline\newline\noindent
\subsection{Statement of contribution.}
With respect to the large literature outlined above, we focus on the specific Problem~(\ref{obj_fcn})--(\ref{constr:init_cond}). By varying $\lambda$, we can change the weight of the minimum-time and minimum-fuel tasks in the objective function. Problem~(\ref{obj_fcn})--(\ref{constr:init_cond}) is based on some modelling simplification (that we will discuss in Appendix~\ref{sec:prob_deriv}). The overall model is simpler than some models used in the above literature. However, under mild assumptions, we can solve Problem~(\ref{obj_fcn})--(\ref{constr:init_cond}) very efficiently, with the guarantee of finding the optimal solution. Namely, our main contributions are the following ones:
\begin{enumerate}
\item We present a convex relaxation of this problem (see, in what follows, Problem~(\ref{obj_fcn_t})--(\ref{constr:forces_f})). We present mild conditions that guarantee that this relaxation is exact (Propositions~\ref{prop:condwmax},~\ref{prop:condwbar}). This property guarantees that we can find the global optimum of Problem~(\ref{obj_fcn})--(\ref{constr:init_cond}) with convex solvers. Moreover, as we will show in the simulation examples (Section~\ref{sec:num_exp}), the solution is very efficient in terms of computational time. 
%  solving efficiently Problem~(\ref{obj_fcn})--(\ref{constr:init_cond})
%w
 % provides the exact solution to the original non-convex speed planning problem, meaning that we can exactly sol%ve a non-convex problem by means of efficient convex solvers. Such results
%are hidden convexity results, where the best known of such results are the one related to the hidden convexity of the trust region problem, i.e., the problem of minimizing a non-convex quadratic function over the unit sphere, and of
%the so called generalized trust region problem, where a quadratic objective function
%is minimized over a region defined by a single, not necessarily convex, quadratic inequality constraint.
%For the latter problem, from a result derived in~\cite{fradkov1979thes} it is possible to prove that the convex Shor relaxation (see~\cite{shor1987quadratic}) is exact.  
\item We present a counterexample (Section~\ref{sec:counter}) showing that, if the previous conditions are not met, the relaxed problem might not be exact.
\end{enumerate}
We present numerical experiments with a focus on the Pareto front of the trade-off between the two objectives (travel time and energy consumption), and the differences between the use of a thermal engine vehicle and an electric vehicle.
\subsubsection{Paper organization}
This work is structured as follows: in Section~\ref{sec:formul}, we introduce a reformulation and convex relaxations of Problem (\ref{obj_fcn})--(\ref{constr:init_cond}). Section~\ref{sec:exact} presents the main theoretical results, i.e., the exactness conditions for the convex relaxation
of the non-convex speed planning problem (Section~\ref{sec:theo}), and a counterexample showing that the relaxation might not be exact when the exactness conditions do not hold (Section~\ref{sec:counter}). Section~\ref{sec:num_exp} shows the numerical experiments. Section~\ref{sec:concl} draws the conclusions and highlights possible future directions. Finally, in the Appendix~\ref{sec:prob_deriv} one can find all the details of the problem derivation.
\section{Convex relaxation of Problem (\ref{obj_fcn})--(\ref{constr:init_cond})}
\label{sec:formul}
We reformulate problem~(\ref{obj_fcn})--(\ref{constr:init_cond}), introducing auxiliary variables $t_i$, with $i \in \{1, \ldots, n-1\}$, as follows:
\begin{align}
& \min_w \sum_{i=1}^{n-1} h \left(\lambda \max\{\eta F_i,F_i\} + t_i\right) & \label{obj_fcn_t} \\
\textrm{such that} & & \nonumber \\
& w_i \leq \wmax_i & \text{for } i \in \{1, \ldots, n\}, \nonumber \\
& |F_i| \leq Mg\mu & \text{for } i \in \{1, \ldots, n-1\}, \nonumber \\
& t_i = \frac{1}{\sqrt{w_i}}& \text{for } i \in \{1, \ldots, n-1\}, \label{constr:t} \\
& t_i \geq \frac{F_i}{\Pmax} & \text{for } i \in \{1, \ldots, n-1\}, \label{constr:max_power_t} \\
& \frac{M}{h} (w_{i+1}-w_i) = - \Gamma w_i + F_i - M g (\sin \alpha_i + c) & \text{for } i \in \{1, \ldots, n-1\}, \nonumber \\
& w_i \geq 0 & \text{for } i \in \{1, \ldots, n\}, \nonumber \\
& w_1 = \winit, & \nonumber
\end{align}
where constraints~\eqref{constr:t} define the new auxiliary variables $t_i$, while~\eqref{obj_fcn_t} and constraints~\eqref{constr:max_power_t} are just a rewriting of objective function~\eqref{obj_fcn} and constraints~\eqref{constr:max_power} in terms of $t_i$ and $F_i$, respectively. In this new formulation non-convexity is moved into the equality constraints~\eqref{constr:t}.
Now, for ease of notation, by setting $f_i := \frac{F_i}{M}$, for $i \in \{1, \ldots, n-1\}$, and $\gamma := \frac{\Gamma}{M}$, we have the following reformulation:
\begin{align}
& \min_w \sum_{i=1}^{n-1} h \left(\lambda M\max\{\eta f_i,f_i\} + t_i\right) & \nonumber \\
\textrm{such that} & & \nonumber \\
& w_i \leq \wmax_i & \text{for } i \in \{1, \ldots, n\}, \nonumber \\
& |f_i| \leq g\mu & \text{for } i \in \{1, \ldots, n-1\}, \nonumber \\
& t_i = \frac{1}{\sqrt{w_i}}& \text{for } i \in \{1, \ldots, n-1\}, \nonumber \\
& t_i \geq \frac{M f_i}{\Pmax} & \text{for } i \in \{1, \ldots, n-1\}, \nonumber \\
& f_i = \frac{1}{h} (w_{i+1}-w_i) + \gamma w_i + g (\sin \alpha_i + c) & \text{for } i \in \{1, \ldots, n-1\}, \nonumber \\
& w_i \geq 0 & \text{for } i \in \{1, \ldots, n\}, \nonumber \\
& w_1 = \winit. & \nonumber
\end{align}
This can be relaxed into the following convex problem by simply turning the equality constraints~\eqref{constr:t} into inequality constraints: 
\begin{align}
& \min_w \sum_{i=1}^{n-1} h \left(\lambda M\max\{\eta f_i,f_i\} + t_i\right) & \nonumber \\
\textrm{such that} & & \nonumber \\
& w_i \leq \wmax_i & \text{for } i \in \{1, \ldots, n\}, \nonumber \\
& |f_i| \leq g\mu& \text{for } i \in \{1, \ldots, n-1\}, \label{constr:max_force_f} \\
& t_i \geq \frac{1}{\sqrt{w_i}}& \text{for } i \in \{1, \ldots, n-1\}, \label{constr:t_relax} \\
  & t_i \geq \frac{M f_i}{\Pmax} & \text{for } i \in \{1, \ldots, n-1\},  \label{constr:max_power_tf1} \\
& f_i = \frac{1}{h} (w_{i+1}-w_i) + \gamma w_i + g (\sin \alpha_i + c) & \text{for } i \in \{1, \ldots, n-1\}, \label{constr:forces_f} \\
& w_i \geq 0 & \text{for } i \in \{1, \ldots, n\}, \nonumber \\
& w_1 = \winit. & \nonumber
\end{align}
Constraints~\eqref{constr:t_relax} are now convex constraints which relax constraints~\eqref{constr:t}.
The previous convex relaxation can be formulated as a Second Order Cone Programming (SOCP) problem by rewriting constraints~\eqref{constr:t_relax} as $1 \leq t_i w_i^{-\frac{1}{2}}$, which are equivalent to 
\[
1 \leq z_i y_i,\quad y_i^2 \leq t_i,\quad z_i^2 \leq t_i w_i,
\]
(see, e.g., ~\cite{alizadehSecondorderConeProgramming2003}).
These are SOCP constraints. In general, any constraint of form $a^2 \leq b c$, with $a,b,c$ positive is equivalent to
\[
\left\| \begin{array}{cc} 2 a \\ b-c\end{array} \right\| \leq b+c.
\]

\section{Main results}
\label{sec:exact} 
In this section, we address the following:
\begin{qstn}
Can we prove that, under suitable assumptions, at optimal solutions of the relaxed problem, constraints~\eqref{constr:t_relax} are always active and, thus, that the relaxed problem is exact?
\end{qstn}
In what follows, we will derive exactness conditions and a counterexample for the case in which such conditions are not met.

\subsection{Exactness conditions}
\label{sec:theo}
First of all, we define the critical (squared) speed $\bar w$ as the solution of
\begin{equation}\label{eq:critical_speed}
\frac{\Pmax}{M\sqrt{\bar w}} = g\mu.
\end{equation}
This quantity will be essential in what follows.
Note that the critical speed is the one at which the maximum traction constraint~\eqref{constr:max_force} and the maximum power consumption constraint~\eqref{constr:max_power} take on the same value with respect to force $F_i$, i.e., the two constraints are both active. Physically speaking, the critical speed is the one at which, using the maximum available power, the tires are at the limit of slipping due to static friction, in the absence of other forces such as rolling resistance and aerodynamic drag.
At optimal solutions of the relaxed problem, at least one of the two constraints (\ref{constr:t_relax}) and (\ref{constr:max_power_tf1}) is active. Moreover, we make the following remark.

\begin{rmrk}
  \label{rem:viol}
Let $(w^*, t^*, f^*)$ be the optimal solution of the relaxed problem. If the convex relaxation is not exact, there exists $i \in \{1, \ldots, n-1\}$ such that constraint~\eqref{constr:t_relax} is not active, and $i$ is the last index for which this happens. This means that
\begin{equation}
\label{ineq:t_violation}
t_i^* = \frac{M f_i^*}{\Pmax} > \frac{1}{\sqrt{w_i^*}}\enspace \wedge\enspace (\forall j\in \{i+1,\ldots,n-1\})\ t_j^*= \frac{1}{\sqrt{w_j^*}}\ .
\end{equation}
\end{rmrk}
We first prove that~\eqref{ineq:t_violation} cannot hold at a speed lower or equal than the critical one.
\begin{proposition}
\label{prop:underbarw}
Let $(w^*, t^*, f^*)$ be the optimal solution of the relaxed problem. Then, for all $i \in \{1, \ldots, n\}$ such that $w_i^* \leq \bar w$, with $\bar w$ defined as in~\eqref{eq:critical_speed}, constraints~\eqref{constr:t_relax} are active.
\end{proposition}

\begin{proof}
Let us assume by contradiction that~\eqref{ineq:t_violation} holds. 
%that for some $i \in \{1, \ldots, n-1\}$, constraint~\eqref{constr:t_relax} is not active, and that $i$ is the last index for which this happens. This means that
%\begin{equation}
%\label{ineq:t_violation}
%t_i^* = \frac{M f_i^*}{\Pmax} > \frac{1}{\sqrt{w_i^*}}\enspace \wedge\enspace (\forall j\in \{i+1,\ldots,n-1\})\ t_j^*= \frac{1}{\sqrt{w_j^*}}\ .
%\end{equation}
We observe that conditions~\eqref{constr:max_force_f} and~\eqref{ineq:t_violation}, together with the definition (\ref{eq:critical_speed}) of the critical speed, imply that
\begin{equation}
\label{ineq:mu_grt}
\frac{\Pmax}{M \sqrt{\bar w}} = g\mu \geq f_i^* > \frac{\Pmax}{M \sqrt{w_i^*}}\ .
\end{equation}
However, if $w_i^* \leq \bar w$, condition~\eqref{ineq:mu_grt} cannot hold, implying that constraint~\eqref{constr:t_relax} is active.
\end{proof}

\begin{proposition}
\label{prop:f_leq_mu}
Let $(w^*, t^*, f^*)$ be a feasible solution of the relaxed problem, and assume that $h$ satisfy
\begin{equation}
\label{ineq:h_suff_cond}
(1 - h\gamma) \bar w - hg (1 + c) > \left(\frac{\Pmax h}{2M(\lambda \gamma \Pmax h + 1 - \lambda)}\right)^{\frac{2}{3}},
\end{equation}
which always holds if the length step $h$ is small enough.
If there exists $i \in \{1, \ldots, n-2\}$ such that~\eqref{ineq:t_violation} holds true and
\begin{equation}
\label{ineq:f_less_mu}
f_{i+1}^* < g\mu,
\end{equation}
then $(w^*, t^*, f^*)$ cannot be the optimal solution of the relaxed problem.
\end{proposition}

\begin{proof}
Let us consider a new solution $(\tilde w, \tilde t, \tilde f)$, which coincides with $(w^*, t^*, f^*)$, except that, for some small $\delta > 0$,
\begin{align*}
\tilde w_{i+1}	& = w_{i+1}^* - \delta \\
\tilde f_i		& = f_i^* - \frac{\delta}{h} \\
\tilde f_{i+1}	& = f_{i+1}^* + \frac{\delta}{h} -\delta\gamma \\
\tilde t_i		& = t_i^* - \frac{M}{\Pmax}\frac{\delta}{h} \\
\tilde t_{i+1}	& = \max\left\{ \frac{M}{\Pmax}\left( f_{i+1}^* + \frac{\delta}{h} -\delta\gamma \right), \frac{1}{\sqrt{\tilde w_{i+1}}} \right\} \approx \\
			& \approx \max\left\{ \frac{M}{\Pmax}\left( f_{i+1}^* + \frac{\delta}{h} -\delta\gamma \right), \frac{1}{\sqrt{w_{i+1}^*}} + \frac{\delta}{2 (w_{i+1}^*)^{\frac{3}{2}}} \right\}.
\end{align*}
Condition~\eqref{ineq:f_less_mu} ensures that we can consider a sufficiently small $\delta$ such that also $\tilde f_{i+1}$ satisfies~\eqref{constr:max_force_f}, that is, $\tilde f_{i+1} \leq g\mu$, so that $(\tilde w, \tilde t, \tilde f)$ is feasible for the relaxed problem.
Now, let us consider the difference between the objective function values at $(\tilde w, \tilde t, \tilde f)$ and at $(w^*, t^*, f^*)$, knowing that $f_i^*, f_{i+1}^*, \tilde f_i, \tilde f_{i+1} \geq 0$,
\begin{gather*}
\lambda M(\tilde f_i + \tilde f_{i+1} - f_i^* - f_{i+1}^*) + (\tilde t_i + \tilde t_{i+1} - t_i^* - t_{i+1}^*) = \\
= - \lambda M \delta \gamma +  \left(- \frac{M}{\Pmax}\frac{\delta}{h} + \max\left\{ \frac{M}{\Pmax}\left( f_{i+1}^* + \frac{\delta}{h} -\delta\gamma \right), \frac{1}{\sqrt{w_{i+1}^*}} + \frac{\delta}{2 (w_{i+1}^*)^{\frac{3}{2}}} \right\} - t_{i+1}^*\right).
\end{gather*}
Recalling that
\[
t_{i+1}^* = \max \left\{ \frac{M f_{i+1}^*}{\Pmax}, \frac{1}{\sqrt{w_{i+1}^*}} \right\},
\]
we have that
\begin{gather*}
\lambda M(\tilde f_i + \tilde f_{i+1} - f_i^* - f_{i+1}^*) + (\tilde t_i + \tilde t_{i+1} - t_i^* - t_{i+1}^*) \leq \\
\leq - \lambda M \delta \gamma + \max\left\{ -\frac{M}{\Pmax}\delta\gamma, - \frac{M}{\Pmax}\frac{\delta}{h} + \frac{\delta}{2 (w_{i+1}^*)^{\frac{3}{2}}} \right\} .
\end{gather*}
Finally, observe that $- \lambda M \delta \gamma - \frac{M}{\Pmax}\delta\gamma < 0$, whilst the other term is negative if
\begin{equation}
\label{ineq:w_exact_cond}
- \lambda M \delta \gamma - \frac{M}{\Pmax}\frac{\delta}{h} + \frac{\delta}{2 (w_{i+1}^*)^{\frac{3}{2}}} < 0
\quad \Longleftrightarrow \quad
w_{i+1}^* > \left(\frac{\Pmax h}{2M(\lambda \gamma \Pmax h + 1)}\right)^{\frac{2}{3}}.
\end{equation}
This implies that, if $w_{i+1}^*$ satisfies~\eqref{ineq:w_exact_cond}, then
\[
\lambda M(\tilde f_i + \tilde f_{i+1} - f_i^* - f_{i+1}^*) + (\tilde t_i + \tilde t_{i+1} - t_i^* - t_{i+1}^*) < 0,
\]
that is, $(\tilde w, \tilde t, \tilde f)$ is another feasible solution with an objective function value strictly lower than that of $(w^*, t^*, f^*)$, meaning that the latter cannot be the optimal one.
\newline\newline\noindent
Note that
\begin{equation}
\label{ineq:wi+1_lb}
f_i^* \geq 0 \ \Longleftrightarrow\
w_{i+1}^* \geq (1 - h\gamma) w_i^* - hg (\sin \alpha_i + c)\ \Longrightarrow\
w_{i+1}^* \geq (1 - h\gamma) \bar w - hg (1 + c),
\end{equation}
where the last inequality follows from Proposition~\ref{prop:underbarw} and observing that $\sin \alpha_i \leq 1$. 
Now, by putting together~\eqref{ineq:w_exact_cond} and~\eqref{ineq:wi+1_lb}, we obtain that a sufficient condition for~\eqref{ineq:w_exact_cond} to hold true is that $h$ must satisfy the inequality~\eqref{ineq:h_suff_cond}.
%\[
%(1 - h\gamma) \bar w - hg (1 + c) > \left(\frac{\Pmax h}{2M(\lambda \gamma \Pmax h + 1 - \lambda)}\right)^{\frac{2}{3}}.
%\]
Note that %, in the previous inequality, which coincides with~\eqref{ineq:h_suff_cond}, 
the left-hand side of~\eqref{ineq:h_suff_cond} monotonically converges to $\bar w$ from below as $h \rightarrow 0$, whilst the right-hand side monotonically converges to $0$ from above as $h \rightarrow 0$.
%, provided that $\lambda < 1$.
Hence, there exists a sufficiently small $\bar h > 0$ such that, for all $h \in (0, \bar h)$, inequality~\eqref{ineq:h_suff_cond} is satisfied.
%Otherwise, if $\lambda = 1$, that is, we are exclusively optimizing with respect to energy consumption and are completely neglecting travel time, then~\eqref{ineq:h_suff_cond} becomes
%\[
%\bar w - h(\gamma \bar w + g (1 + c)) > \frac{1}{\left(2M\gamma\right)^{\frac{2}{3}}}\ \Longleftrightarrow\
%h < \frac{\left(2\Gamma\right)^{\frac{2}{3}}\bar w - 1}{\left(2\Gamma\right)^{\frac{2}{3}}(\gamma \bar w + g (1 + c))},
%\]
%where the right hand-side of the last inequality is significantly greater than $0$ in practical applications, as we will see in the numerical experiments.
\end{proof}

\begin{proposition}
Let $(w^*, t^*, f^*)$ be a feasible solution of the relaxed problem and assume that condition~\eqref{ineq:t_violation} holds for $i = n - 1$, then such solution is not the optimal one.
\end{proposition}

\begin{proof}
Similarly to the proof of Proposition~\ref{prop:f_leq_mu},
we can consider a new solution $(\tilde w, \tilde t, \tilde f)$ feasible for the relaxed problem, identical to $(w^*, t^*, f^*)$ except that, for some small $\delta > 0$, $\tilde w_n = w_n^* - \delta$.
As a consequence, $\tilde f_{n-1} < f_{n-1}^*$ and $\tilde t_{n-1} < t_{n-1}^*$, making $(\tilde w, \tilde t, \tilde f)$ another feasible solution with a smaller objective function value, meaning that $(w^*, t^*, f^*)$ cannot be the optimal one.
\end{proof}

%\begin{align*}
%h	&= \frac{w_i - w_{i+1}}{\gamma w_i + g (\sin \alpha_i + c) - f_i}\\
%	&\frac{w_i - w_{i+1}}{\gamma w_i + g (1 + c) - f_i}
%	\geq \frac{\bar w - w_{i+1}}{\gamma \max\{\wmax, \bar w\} + g (1 + c)}
%\end{align*}

Observe that Proposition~\ref{prop:f_leq_mu} states that, for the optimal solution of the relaxed problem, if~\eqref{ineq:t_violation} holds true, then~\eqref{ineq:f_less_mu} does not, that is, $f_{i+1}^* = g\mu$.
Then, if~\eqref{ineq:t_violation} holds true, we have that
\begin{equation}
\label{ineq:mu_leq}
g\mu\leq \frac{\Pmax t_{i+1}^*}{M} = \frac{\Pmax}{M \sqrt{w_{i+1}^*}}.
\end{equation}
As a result, by putting together inequalities~\eqref{ineq:mu_leq} and~\eqref{ineq:mu_grt}, we obtain
\begin{equation}
\label{ineq:w_decr}
\frac{\Pmax}{M \sqrt{w_i^*}} < g\mu \leq \frac{\Pmax}{M \sqrt{w_{i+1}^*}}
\quad \Longleftrightarrow \quad
\sqrt{w_i^*} > \sqrt{w_{i+1}^*},
\end{equation}
meaning that, from $i$ to $i+1$, speed must be strictly decreasing.
At this point, if we consider constraints~\eqref{constr:t_relax}--\eqref{constr:forces_f} at index $i$, we obtain, when~\eqref{ineq:t_violation} holds true, that:
\begin{equation}
\label{ineq:not_exact}
\frac{\Pmax}{M\sqrt{w_i^*}} < \frac{1}{h} (w_{i+1}^* - w_i^*) + \gamma w_i^* + g (\sin \alpha_i + c)
\end{equation}
Then, 
any condition under which (\ref{ineq:not_exact}) cannot hold is also a condition under which the relaxed convex problem provides an optimal solution of the original non-convex problem~\eqref{obj_fcn}--\eqref{constr:init_cond}.
We provide two such conditions. The first one, reported in Proposition~\ref{prop:condwmax}, refers to the maximum speed profile, while the second, reported in Proposition~\ref{prop:condwbar}, refers to the critical speed.
\begin{proposition}
\label{prop:condwmax}
Let us assume that (\ref{ineq:h_suff_cond}) holds.
If the maximum squared speed profile $\wmax$ satisfies
\begin{equation}
\label{ineq:exact_cond}
(\forall i \in \{1, \ldots, n\})\ \frac{\Pmax}{M\sqrt{\wmax_i}} \geq \gamma \wmax_i + g (\sin \alpha_i + c),
\end{equation}
then, the optimal solution of the relaxed problem is feasible and, thus, optimal for Problem~\eqref{obj_fcn}--\eqref{constr:init_cond}.
\end{proposition}

\begin{proof}
Since at any feasible solution to both the original problem~\eqref{obj_fcn}--\eqref{constr:init_cond} and the relaxed problem it holds that $w_i \leq \wmax_i$, for $i \in \{1, \ldots, n\}$, then:
\[
\frac{\Pmax}{M\sqrt{w_i}} \geq \frac{\Pmax}{M\sqrt{\wmax_i}}.
\]
Moreover, note that if $w_{i+1} < w_i$, then $\frac{1}{h} (w_{i+1} - w_i) < 0$, so if condition~\eqref{ineq:exact_cond} holds true, we obtain that, for any feasible solution of the relaxed problem, for all $i \in \{1, \ldots, n-1\}$ such that $w_{i+1} < w_i$, it holds that,
\[
\frac{\Pmax}{M\sqrt{w_i}} \geq \frac{\Pmax}{M\sqrt{\wmax_i}} \geq \gamma \wmax_i + g (\sin \alpha_i + c) \geq \gamma w_i + g (\sin \alpha_i + c) \geq \frac{1}{h} (w_{i+1} - w_i) + \gamma w_i + g (\sin \alpha_i + c),
\]
meaning that condition~\eqref{ineq:not_exact}, which is implied by~\eqref{ineq:t_violation} and~\eqref{ineq:mu_leq}, cannot hold.
\end{proof}
\begin{proposition}
\label{prop:condwbar}
Let us assume that (\ref{ineq:h_suff_cond}) holds.
If:
\begin{equation}
\label{eq:condbarw}
\min_{i\in I} \left[\frac{\Pmax}{M}\sqrt{\frac{1 - h\gamma}{\bar w + hg(\sin \alpha_i + c)}}-\frac{\gamma}{1-h\gamma}(\bar w + hg(\sin \alpha_i + c))- g (\sin \alpha_i + c)\right]\geq 0,
\end{equation}
where
$$
I=\{i\in \{1,\ldots,n\}\ :\ \wmax_i>\bar{w}\},
$$
then, the optimal solution of the relaxed problem is feasible and, thus, optimal for Problem~\eqref{obj_fcn}--\eqref{constr:init_cond}.
\end{proposition}
\begin{proof}
First, we observe that we can restrict the attention to $i\in I$ since according to Proposition~\ref{prop:underbarw}, for $i\not\in I$
(\ref{ineq:t_violation}) cannot hold. Now, let $(w^*, t^*, f^*)$ be an optimal solution of the relaxed problem and assume by contradiction that condition~\eqref{ineq:t_violation} holds for some $i\in I$.
If~\eqref{ineq:t_violation} holds true, we have previously observed that (\ref{ineq:mu_leq}) holds. Recalling the definition of $\bar w$, we also have that
$$
\frac{\Pmax}{M \sqrt{w_{i+1}^*}}\geq g\mu= \frac{\Pmax}{M \sqrt{\bar w}},
$$
so that $w_{i+1}^*\leq \bar w$.
Next, we observe that
$$
f_i^* = \frac{\Pmax t_i^*}{M} = \frac{1}{h} (w_{i+1}^* - w_i^*) + \gamma w_i^* + g (\sin \alpha_i + c) > \frac{\Pmax }{M \sqrt{w_i^*}} >0,
$$
implies
$$
w_i^* < \frac{1}{(1 - h\gamma)}\left[w_{i+1}^* + hg(\sin \alpha_i + c)\right]
	\leq \frac{1}{(1 - h\gamma)}\left[\bar w + hg(\sin \alpha_i + c)\right].
$$
Then:
\begin{equation}
\label{eq:interm}
\frac{\Pmax}{M\sqrt{w_i^*}}- \gamma w_i^* > \frac{\Pmax}{M}\sqrt{\frac{1 - h\gamma}{\bar w + hg(\sin \alpha_i + c)}}-\frac{\gamma}{1-h\gamma}(\bar w + hg(\sin \alpha_i + c)).
\end{equation}
Recalling that~\eqref{ineq:t_violation} implies (\ref{ineq:not_exact}), we have that:
$$
 \frac{1}{h} (w_{i+1}^* - w_i^*) >\frac{\Pmax}{M\sqrt{w_i^*}} - \gamma w_i^* - g (\sin \alpha_i + c),
$$
where, in view of (\ref{eq:interm}), the right-hand side is not negative for each $i\in I$ under condition (\ref{eq:condbarw}). 
Thus, we are led to a contradiction since, as previously observed, we should have that $w_{i+1}^*<w_i^*$ when~\eqref{ineq:t_violation} holds true.
\end{proof}
Note that for $h\rightarrow 0$,  (\ref{ineq:h_suff_cond}) holds, while condition (\ref{eq:condbarw}) reduces to:
$$
\min_{i\in I} \left[\frac{\Pmax}{M\sqrt{\bar w}}-\gamma\bar w- g (\sin \alpha_i + c)\right]\geq 0,
$$
so that we can make the following remark, stating the physical interpretation of condition (\ref{eq:condbarw}).
\begin{rmrk}
For $h\rightarrow 0$ condition (\ref{eq:condbarw}) is fulfilled if, at the critical speed $\bar w$, the vehicle has enough power to accelerate or, at least, not to decelerate, despite drag, positive slope, and tires rolling resistance, along all portions of the path where the maximum allowed speed
exceeds the critical speed.
\end{rmrk} 
%So, if condition~\eqref{ineq:exact_cond} holds, it ensures that an optimal solution to the relaxed problem is feasible for the original problem~\eqref{obj_fcn}--\eqref{constr:init_cond}, since condition~\eqref{ineq:t_violation} is ruled out.
%Observe that, from the practical application point of view, condition~\eqref{ineq:exact_cond} states that, along the whole path, the vehicle is always able to accelerate despite drag, positive gradients, and tires rolling resistance.
%
%\begin{thrm}
%{\color{red} Condizione per $\min\{\wmax, \bar w\}$.}
%\end{thrm}
Let us now show a counterexample in which the relaxed problem is not exact and, indeed, conditions~(\ref{ineq:exact_cond})~\eqref{eq:condbarw} do not hold.

\subsection{Counterexample}
\label{sec:counter}

We present a counterexample in which conditions~(\ref{ineq:exact_cond})~\eqref{eq:condbarw} are not satisfied, and the relaxed problem is not exact. In particular, we look for an instance of Problem~(\ref{obj_fcn})--(\ref{constr:init_cond}) that is not feasible, but it is feasible if we remove the maximum power constraints~(\ref{constr:max_power}). For such an instance, relaxation~(\ref{obj_fcn_t})--~(\ref{constr:forces_f}) is admissible, since the power constraints has been relaxed. Hence, the relaxation cannot be exact.
In particular, we consider a scenario with a steep slope. Due to the constraint on maximum force~(\ref{constr:max_force}), the vehicle is unable to maintain its speed on the slope, since the maximum applicable force is unable to compensate the force due to gravity. Hence, the vehicle has to accelerate before the slope, to reach a sufficiently large velocity. If the power constraint~(\ref{constr:max_power}) is enforced, the vehicle cannot sufficiently accelerate, and Problem~(\ref{obj_fcn})--(\ref{constr:init_cond}) is not feasible. Instead, without the maximum power constraint, the problem is feasible. The relaxed problem is feasible, since solutions can violate the maximum power constraint.
Note that we managed to build the counterexample by considering a vehicle with a very small maximum power, compared to a standard road vehicle. Indeed, the cars on the market have sufficient power to drive uphill at relevant slopes.
\begin{figure}[h]
    \centering
    \includegraphics[width=0.6\textwidth, trim=50 13 50 30, clip]{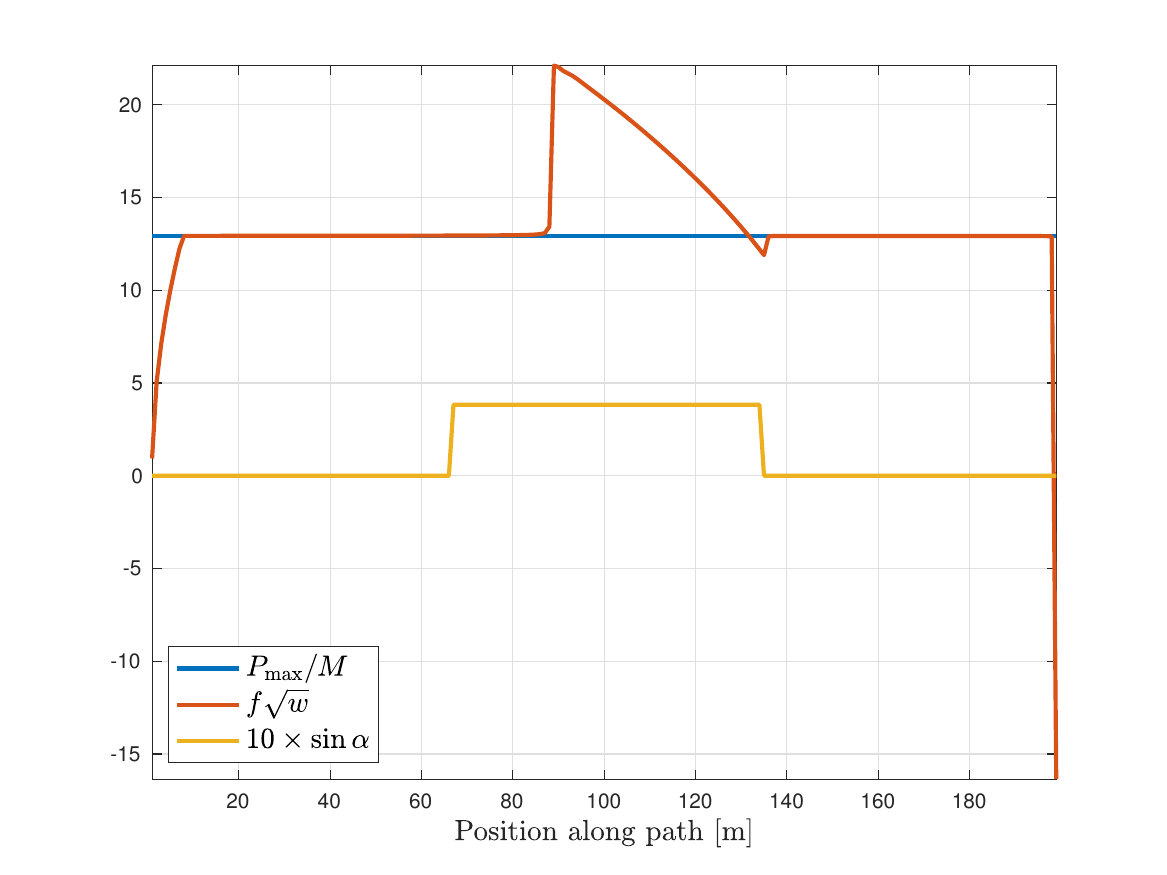}
    \caption{Violation of maximum power constraint.}
    \label{fig:Pmax_viol}
\end{figure}
\begin{figure}[h]
    \centering
    \includegraphics[width=0.6\textwidth, trim=50 13 45 30, clip]{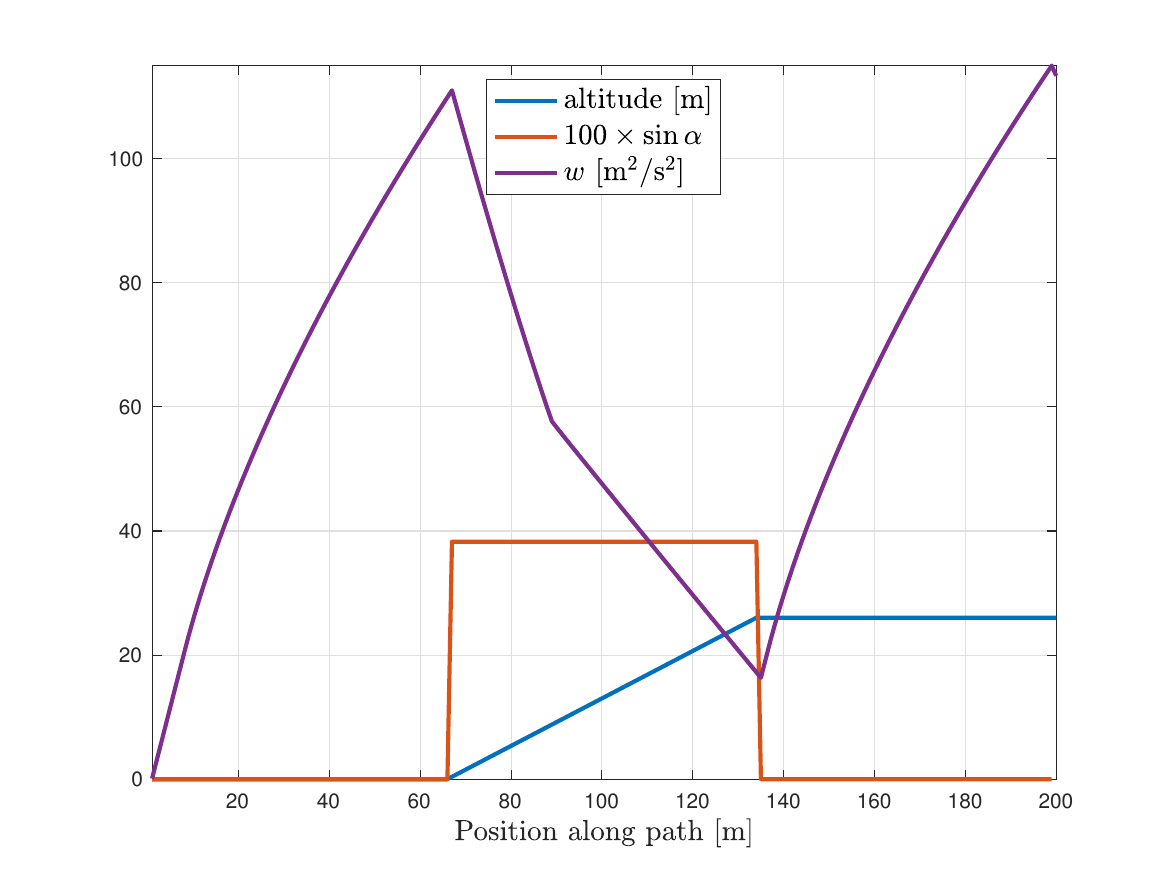}
    \caption{Squared speed profile.}
    \label{fig:w_profile}
\end{figure}

The vehicle we consider is a Fiat 500 with thermal engine, for which we reduced its maximum power to \SI{12500}{\watt} ($\approx 17$ metric hp), whilst its actual maximum power is $69$ metric hp.
All the other vehicle specifications required for the problem are presented in Section~\ref{sec:num_exp}.
We also considered $n = 200$, $h = 1$ (so the path is \SI{200}{\meter} long), $\lambda = 0$ (so we only optimize travel time), $\mu = 0.3$ (representing the static friction coefficient between tires and wet asphalt), and $w_1 = \SI{0.1}{\meter^2\per\second^2}$.
The path we consider is straight and is divided in three sections of length $\approx$ \SI{66.7}{\meter} each, the initial and the final one are both flat, whilst the central section has a constant $22.5^\circ$ (approximately $41.5\%$) incline uphill.
The considered maximum squared speed profile is identically equal to \SI{1975}{\meter^2\per\second^2} (\SI{160}{\kilo\meter\per\hour}) along the whole path.
Note that under these conditions, inequality~\eqref{ineq:h_suff_cond} is satisfied, whilst both conditions~\eqref{ineq:exact_cond} and~\eqref{eq:condbarw} are not.
In particular, the fact that~\eqref{eq:condbarw} is not satisfied can be physically interpreted as follows.
Given the limited maximum power and the extreme incline, in the relaxed problem the vehicle is forced to violate the maximum power constraint in order to be able to make it to the end of the incline and not to slide back down the slope.
One can see the constraint violation in Figure~\ref{fig:Pmax_viol}, in the last portion of the first section of the path, right before the beginning of the incline, and during the entire section of the incline.
This violation is necessary as one can see in Figure~\ref{fig:w_profile}, where the squared speed profile approaches zero at the end of the incline, indeed, its minimum value is \SI{16.35}{\meter^2\per\second^2}, which corresponds to a speed of \SI{14.55}{\kilo\meter\per\hour}.
Had not the vehicle exceeded the maximum power constraint, it would not have made it to the end of the slope.
In both figures, the sine of the path grade is always displayed for reference.
As we just saw, the reason why in the previous counterexample the relaxed problem is not exact is that the original problem~\eqref{obj_fcn}--\eqref{constr:init_cond} is not feasible due to the maximum power constraint. This leads us to the following open question.

\begin{qstn}
Does the feasibility of Problem~\eqref{obj_fcn}--\eqref{constr:init_cond} imply that the relaxed problem is exact?
\end{qstn}
Up to now we have not been able to find an instance showing that the answer to the question above is negative.

\section{Numerical experiments}
\label{sec:num_exp}
In this section we present some numerical experiments that show the Pareto front of the two objective functions and the solution times for different values of parameter $\lambda$ for two types of vehicle engines and also the solution times for different values of $n$ for the same two types of vehicle engines.
For the simulations we considered the two following vehicles: a Fiat 500 with thermal engine and a Fiat 500e full electric.
All the vehicle specification needed for running the numerical tests are reported in Table~\ref{tab:vehicle_spec}.
Note that for the Fiat 500 with thermal engine we derived its mass considering the vehicle with a full tank. 
Moreover, we considered $\mu = 0.7$, representing the static friction coefficient between tires and dry asphalt and $h = 3$, so that~\eqref{ineq:h_suff_cond} is satisfied.
Finally, for solving the instances of the problem in all the carried out numerical experiments, we used the optimization solver MOSEK (version 10.1.31)~\cite{mosek} on a hardware with \SI{16}{\giga\byte} of RAM.
\begin{table}[h]
    \centering
    \begin{tabular}{|c|c|c|}
        \hline
        				& \textbf{Fiat 500}					& \textbf{Fiat 500e} \\
       \hline
        $M$			& \SI{967}{\kilogram}					& \SI{1365}{\kilogram} \\
        $\Pmax$		& 69\ metric hp $\approx$ \SI{50750}{\watt}	& 118\ metric hp $\approx$ \SI{87000}{\watt} \\
        $\eta$		& $0\%$							& $70\%$ \\
        $c$			& 0.007							& 0.007 \\
        $\vmax$		&  \SI{160}{\kilo\meter\per\hour}		& \SI{150}{\kilo\meter\per\hour} \\
        $\wmax$		& \SI{1975}{\meter^2\per\second^2}		& \SI{1736}{\meter^2\per\second^2} \\
        $\Gamma$	& \SI{0.406}{\kilogram\per\meter}		& \SI{0.399}{\kilogram\per\meter} \\
        $\gamma$		& \SI{0.4199e-3}{\per\meter}			& \SI{0.2923e-3}{\per\meter} \\
        \hline
    \end{tabular}
    \caption{Vehicles specifications.}
    \label{tab:vehicle_spec}
\end{table}

\begin{figure}[h]
    \centering
    \includegraphics[width=0.6\textwidth, trim=50 12 45 27, clip]{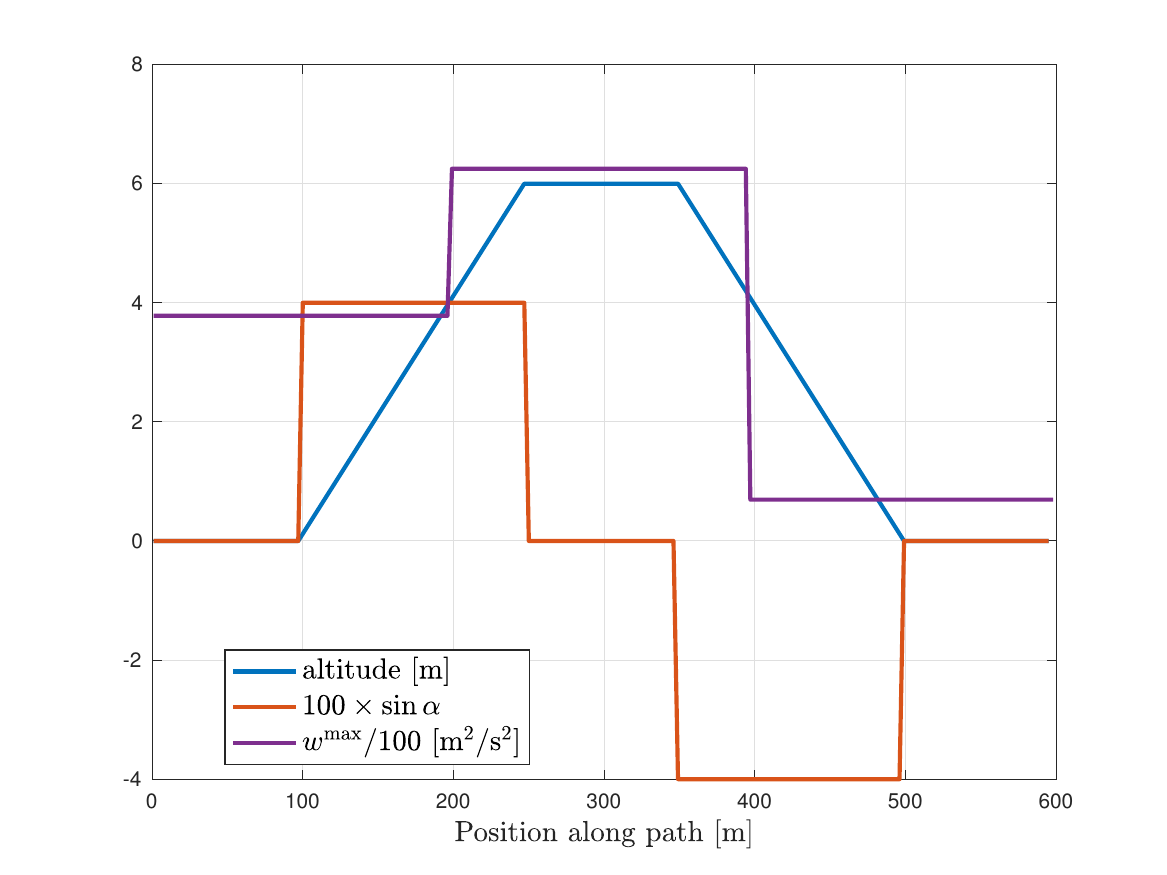}
    \caption{Incline and maximum squared speed profiles of the considered path.}
    \label{fig:altitude}
\end{figure}

For the first set experiments we considered a path of length \SI{600}{\meter} with the following altitude profile: the first \SI{100}{\meter} are flat, followed by a constant $4\%$ positive incline for \SI{150}{\meter} reaching \SI{6}{\meter} of altitude, then the path is flat again for another \SI{100}{\meter}, followed by a constant $-4\%$ negative incline for \SI{150}{\meter}, leading the altitude back to \SI{0}{\meter}, and, finally, one last flat section of \SI{100}{\meter}.
The altitude profile is shown in Figure~\ref{fig:altitude}, together with the sine of the slope angle along the path scaled up by a factor of $100$.
The same figure also displays the maximum squared speed profile $\vmax$, scaled down by a factor of 100, which is the same for both the vehicle with thermal engine and the full electric one, and is structured as follows:
it is divided in three sections of \SI{200}{\meter} each with constant speed, the first one is associated to a maximum speed of \SI{70}{\kilo\meter\per\hour} (\SI{378.1}{\meter^2\per\second^2}), the second to \SI{90}{\kilo\meter\per\hour} (\SI{625.0}{\meter^2\per\second^2}), and the last one to \SI{30}{\kilo\meter\per\hour} (\SI{69.4}{\meter^2\per\second^2}).

\begin{figure}[h]
    \centering
    \includegraphics[width=0.6\textwidth, trim=33 12 50 27, clip]{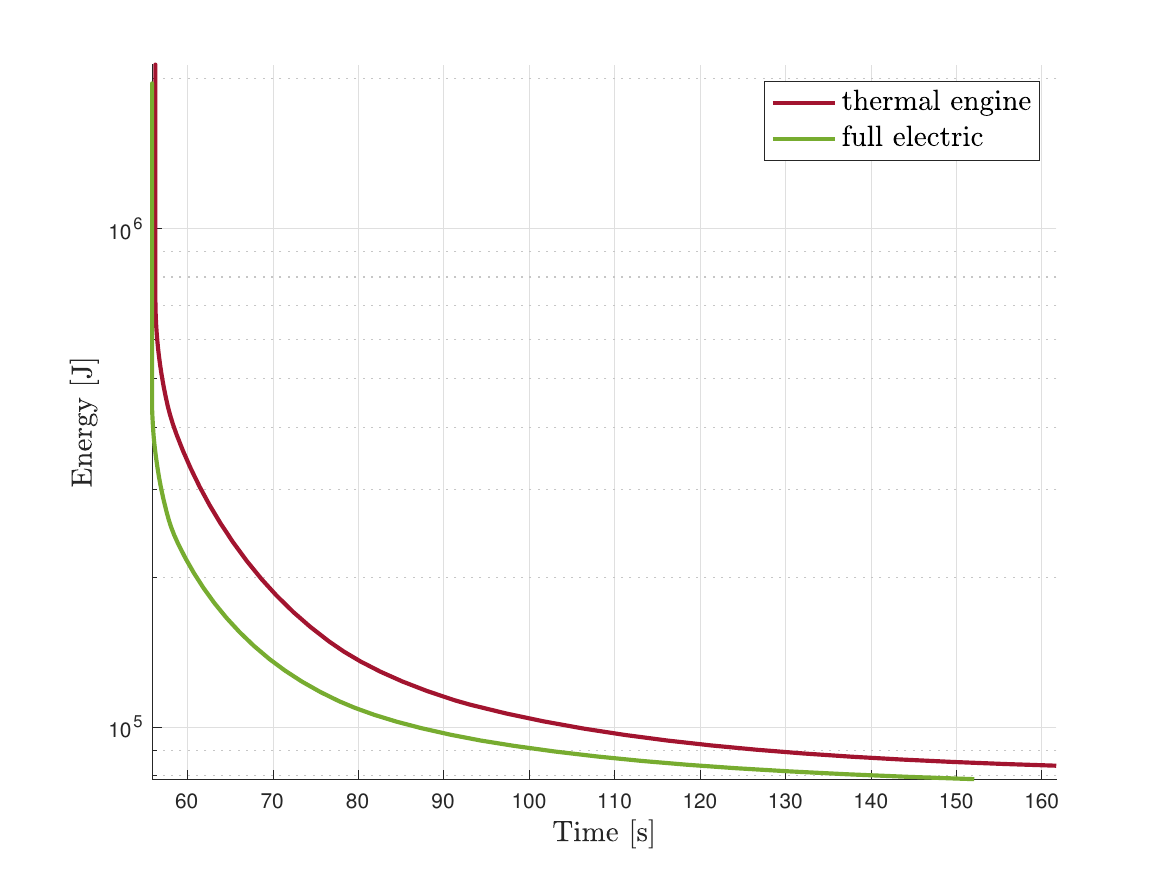}
    \caption{Pareto fronts for the two considered vehicles.}
    \label{fig:T_vs_E}
\end{figure}

\begin{figure}[h]
    \centering
    \includegraphics[width=0.6\textwidth, trim=50 13 45 26, clip]{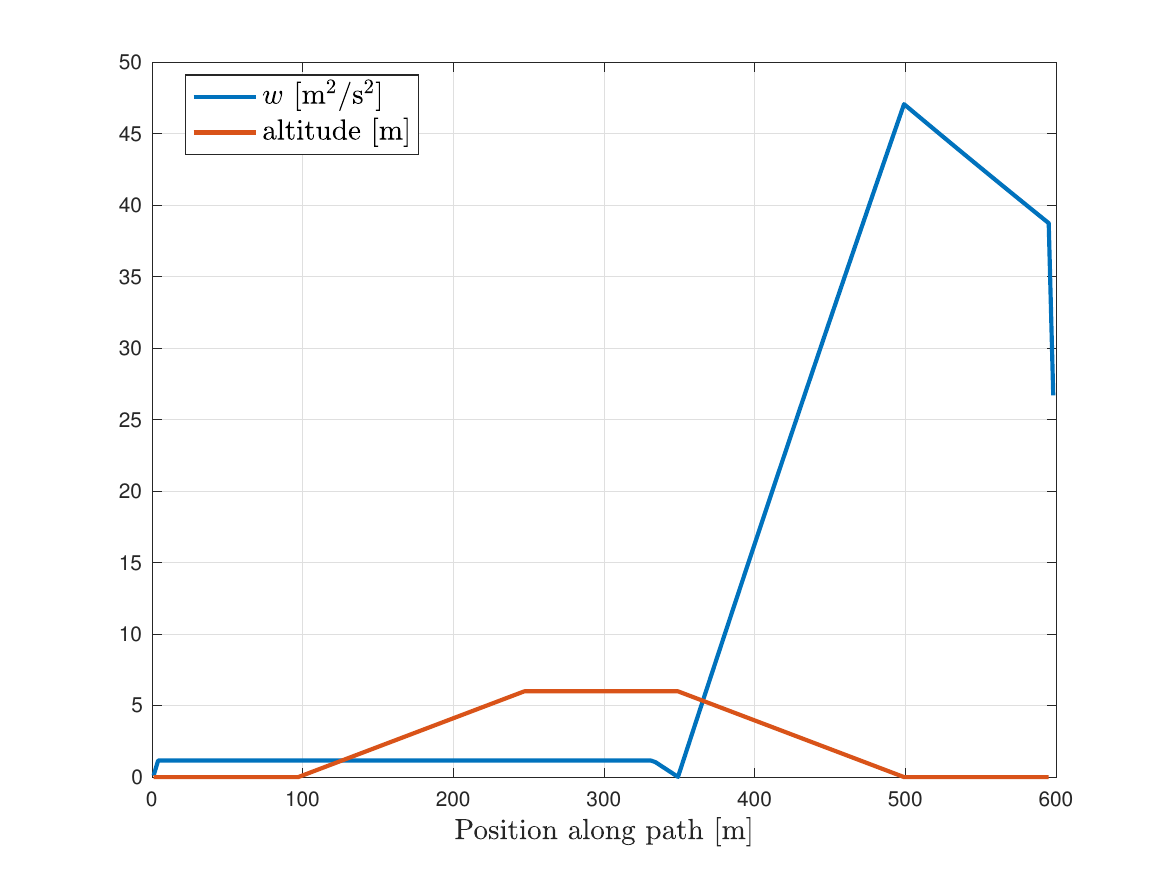}
    \caption{Degenerate squared speed profile for a Fiat 500 with $\lambda = 0.99$.}
    \label{fig:w_dgn}
\end{figure}

For generating the Pareto front of the two objectives we set $n = 200$ and solved the problem at hand for the two considered vehicles for $100$ values of $\lambda$ ranging in $\{0\}$ and from $10^{-7}$ to $10^{-2}$ logarithmically.
As we can see from Figure~\ref{fig:T_vs_E}, the Pareto front associate to the full electric vehicle lies below that of the thermal engine one showing how the Fiat 500e consumes less energy with respect to the Fiat 500 for all values of $\lambda$.
Also note that, in both cases, when we exclusively optimize the travel time, the energy consumption is maximized, whilst on the contrary, as the weight attributed to energy consumption increases, the travel times becomes larger.
However, note that as the value of $\lambda$ approaches $1$, in order to obtain a small decrease in energy consumption, travel times are heavily impacted.
Indeed, if we consider, for instance, $\lambda = 10^{-1}$ we obtain the illustrative squared speed profile shown in Figure~\ref{fig:w_dgn}, for the Fiat 500, in which the average speed in the first half of the path is \SI{3.8}{\kilo\meter\per\hour}, and in the second half, the speed increases up to \SI{24.7}{\kilo\meter\per\hour} only because, in that path section, the vehicle is going down the slope increasing its speed without consuming energy.
This is also the reason why, in the previous experiments, we did not consider values of $\lambda$ greater than $10^{-2}$.
Finally, note that all the previous solutions to the relaxed problems are exact, with a mean maximum violation of equality~\eqref{constr:t} of $8.0\times10^{-8}$ and a maximum one of $6.9\times10^{-7}$.

\begin{figure}[h]
    \centering
    \includegraphics[width=0.6\textwidth, trim=13 25 50 30, clip]{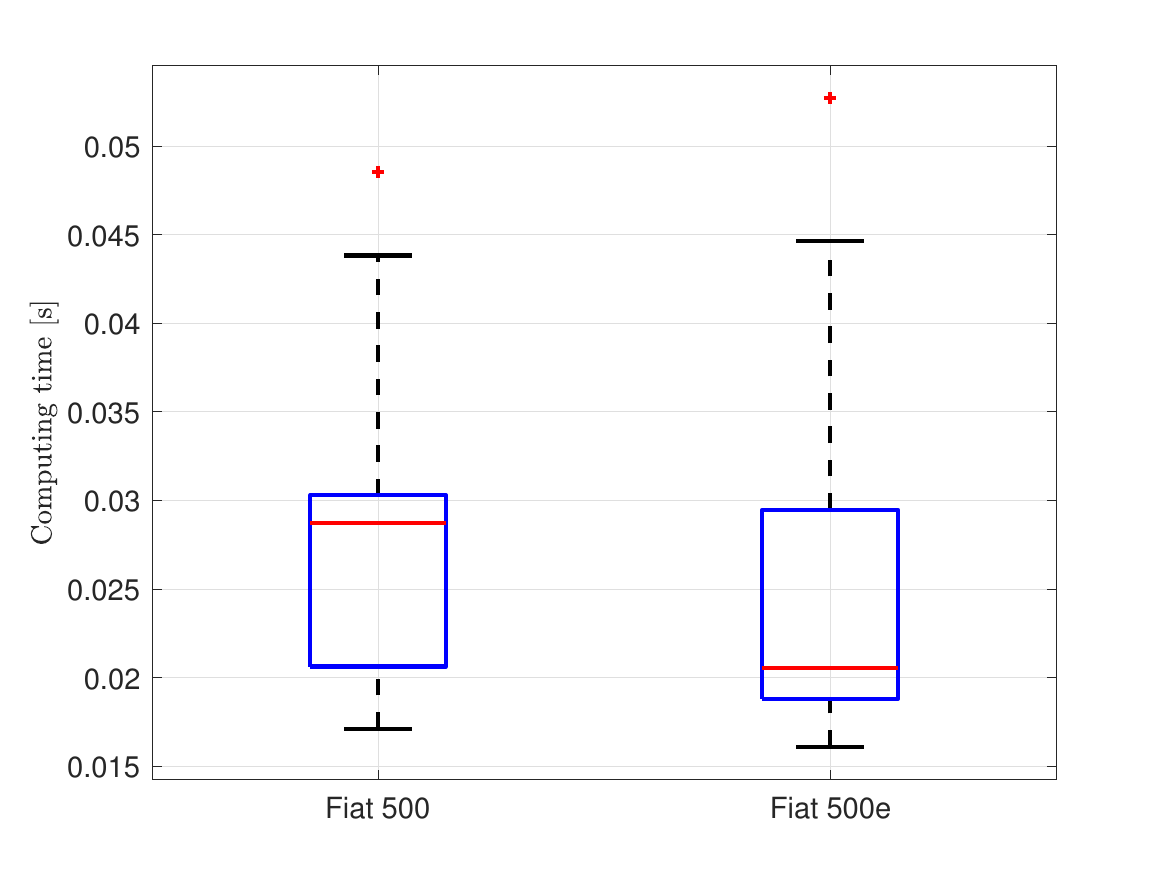}
    \caption{Box-and-whisker plot of the computing times for the two vehicle types.}
    \label{fig:box_plot}
\end{figure}

We also run a set of tests for evaluating the time performances of MOSEK in solving random instances of the problem at hand.
We generated a set of 100 problems for each vehicle type with a smooth random incline profile varying between $5\%$ and $-5\%$, a maximum speed profile divided in three sections of \SI{200}{\meter} each, with constant values for each section randomly taken from $\{30, 50, 70, 90, 110, 130\}$\unit{\kilo\meter\per\hour} $\cup\ \{\vmax\}$, and parameter $\lambda$ that randomly takes value in a set of 100 logarithmically spaced samples in $[10^{-7}, 10^{-2}]$ and $\{0\}$.
The average solver time for solving the previous 100 instances of the problem for the Fiat 500 with thermal engine and for the Fiat 500e full electric are \SI{0.0278}{\second} and \SI{0.0233}{\second}, respectively.
Figure~\ref{fig:box_plot} also shows the box--and-whisker plot of the computing times for solving the previous instances for each vehicle type.

\begin{figure}[h]
    \centering
    \includegraphics[width=0.6\textwidth, trim=30 15 42 27, clip]{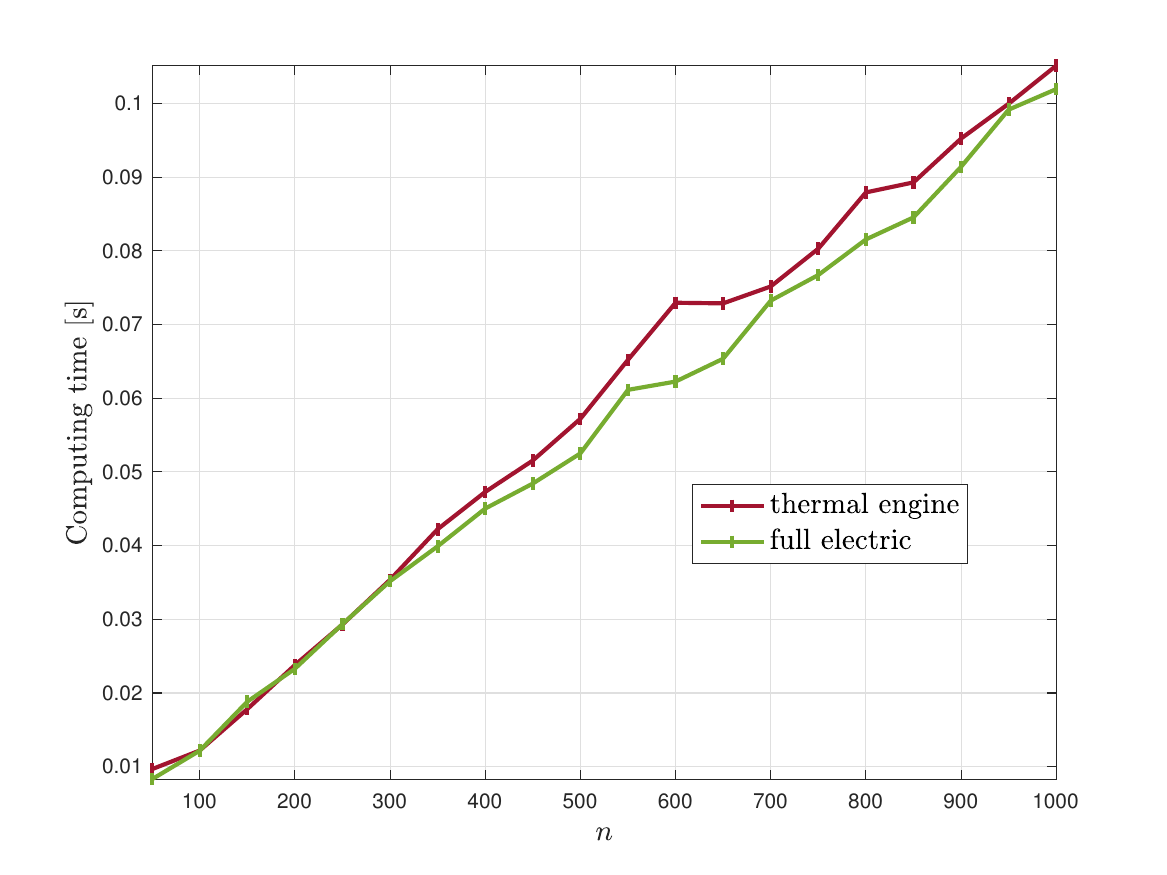}
    \caption{Computing times growth for different values of $n$.}
    \label{fig:n_profile}
\end{figure}

Finally, we also run a set of experiments in order to gain insights into the growth of computing times for increasing values of $n$.
We considered $20$ linearly spaced values of $n$ from $50$ to $1000$.
The altitude and maximum squared speed profiles for these tests are the scaled versions, adapted to $n$ samples, of the profiles shown in Figure~\ref{fig:altitude}.
As we can see from Figure~\ref{fig:n_profile}, the solver times with respect to $n$ exhibit a linear growth for both vehicle engines, with instances of the problem for $n = 1000$ solved in $\approx$ \SI{0.1}{\second}.

\section{Conclusions and future developments}
\label{sec:concl} 
In this paper we addressed the problem of computing the speed profile of a vehicle travelling along an assigned path, subject to suitable constraints (in particular, maximum speed, traction or braking force, and power consumption constraints). A bi-objective optimization approach has been considered since we aim at minimizing at the same time the two conflicting objectives of travel time and energy consumption.
After presenting the resulting mathematical model, which turns out to be a non--convex optimization problem, we introduced a convex relaxation of such problem (in particular, a SOCP relaxation), and we have shown that,
under mild conditions, such relaxation is exact, i.e., it shares the same optimal value of the non--convex problem and an optimal solution of the non--convex problem can be immediately derived from an optimal solution of the relaxation.
Through some numerical experiments, we confirmed that existing solvers for SOCP problems allow solving very efficiently the SOCP relaxation. In the numerical experiments, we have also derived the Pareto front of the trade-off between travel time and energy consumption. We compared the results for two distinct engines, a thermal engine and a full electric engine, where in the latter energy recovery is possible when braking.
\newline\newline\noindent
As possible future developments, we will consider more realistic models of electric vehicles, with a more precise description of regenerative braking and force constraints. We will also explore the possibility of optimizing speed and route choice at the same time, along the lines of our previous work~\cite{9767776}.

\bibliographystyle{plain} 
\bibliography{biblio,VelPlan,VelPlan2}
\appendix
\section{Problem derivation}
\label{sec:prob_deriv}

\begin{figure}
\begin{center}
  \begin{tikzpicture}
    \begin{scope} [rotate=30]
    \shade[top color=red, bottom color=white, shading angle={135}]
        [draw=black,fill=red!20,rounded corners=1.2ex,very thick] (1.5,.5) -- ++(0,1.3) --  ++(4.5,0) -- ++(1.5,-0.5) -- ++(0,-0.8) -- (1.5,.5) -- cycle;
    \draw[very thick, rounded corners=0.5ex,fill=black!20!blue!20!white,thick]  (2.5,1.8) -- ++(0.5,0.7) -- ++(1.8,0) -- ++(1,-0.7) -- (2.5,1.8);
    \draw[thick]  (3.5,1.8) -- (3.5,2.5);
    \draw[draw=black,fill=gray!50,thick] (2.75,.5) circle (.6);
    \draw[draw=black,fill=gray!50,thick] (5.5,.5) circle (.6);
    \draw[draw=black,fill=gray!80,semithick] (2.75,.5) circle (.5);
    \draw[draw=black,fill=gray!80,semithick] (5.5,.5) circle (.5);
    
    \draw[draw=green, thick]  (-0.5,-0.1) coordinate (A1) -- (10,-0.1) coordinate(A2);
    %\draw (8.5,-0.1) node {$x(t)$};
    \draw [fill](4.3,1.3) node (c) {};  
    \node at (2.75,-0.1) (bw) {};
    \node at (5.5,-0.1) (fw) {};
    \draw [-latex,thick] (bw) -- ++(1.0,0) node [above]{$T_r$};
    \draw [-latex,thick] (bw) -- ++(-1.0,0) node [above]{$R_r$};
    \draw [-latex,thick] (fw) -- ++(1.0,0) node [above]{$T_f$};
    \draw [-latex,thick] (fw) -- ++(-1.0,0) node [above]{$R_f$};
    \draw [-latex,thick] (bw) -- ++(0,1.5) node [left]{$W_r$};
    \draw [-latex,thick] (fw) -- ++(0,1.5) node [left]{$W_f$};
    \draw [fill](4.3,1.3) node (c) {} circle (0.1);

%    \draw [dashed] let \p1=(bw) in ($(bw)+(0,0.5)$) -- (\x1,-0.7) node (bq) {} -- (\x1,-0.3);
%  \draw [dashed] let \p1=(fw) in ($(fw)+(0,0.5)$) -- (\x1,-0.7) node (fq) {} (\x1,-0.3);
%  \draw [dashed] let \p1=(c) in (c) -- (\x1,-0.7) node (cq) {} (\x1,-0.3);
  
 % \draw [dashed] let \p1=(c) in (c) -- (-0.1,\y1) node (ch) {} --(-0.2,\y1);

  \node at (9,1) (a) {};
  \draw [-latex] (a)-- ++(-1,0) node [midway,above] {$F_a$};
\end{scope}

\draw (A1)-- ++(10,0);
\draw[thick,blue,-latex] ([shift=(0:10)]A1) arc (0:30:10) node [midway,left] {$\alpha(s)$};
\draw [-latex,thick] (c) -- ++(0,-0.7) node [left,midway] {$M g$};

 \end{tikzpicture}
\end{center}
\caption{Half-car model.}
  \label{fig_car} 
\end{figure}
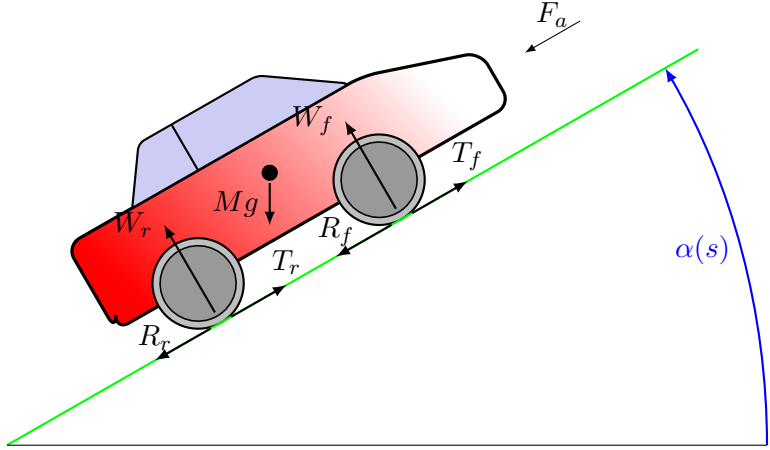

Consider a simplified model of an half-car of mass $M$ and speed $v$, travelling on a road with a variable grade (see Figure~\ref{fig_car}).
The slope angle $\alpha(s)$ is a function of the arc-lenght position.
Here, $F_a(t)$ is the aerodynamic drag force, $T_f$ and $T_r$ are the traction (if positive) or braking (if negative) forces of the front and rear wheels. Forces $R_f$ and $R_r$ are the rolling resistances of the front and rear wheels. We assume that the component of the car acceleration orthogonal to the road is negligible, so that we can write the following dynamic equations, obtained from the balance of forces in the direction parallel and orthogonal to the road:
\begin{equation}
\label{eqn_dyn_1}
  \begin{array}{ll}
    M \dot v(t)= T_f(t)+T_r(t)-R_f(t)-R_r(t)-F_a(t) - Mg \sin \alpha(s(t))\\ [8pt]
0=W_f(t)+W_r(t)-Mg \cos \alpha(s(t)).
\end{array}    
\end{equation}

From the second equation the overall load is $W_f(t)+W_r(t)=Mg \cos \alpha(s(t))$.
To prevent wheels slipping, traction (or braking) forces must satisfy constraints
\[
|T_f| \leq \mu W_f, |T_r| \leq \mu W_r,
\]
assuming that the front and rear axes have the same friction coefficient $\mu$.
From the second of~(\ref{eqn_dyn_1}), we obtain
\begin{equation}
  \label{eqn_no_slip}
|T_f(t)+T_r(t)| \leq \mu M g \cos \alpha(s(t)).
\end{equation}

We assume that the car has four-wheel drive, and that drive and braking forces are optimally balanced. In this case, condition~(\ref{eqn_no_slip}) is sufficient to prevent slipping on front and rear wheels. We assume that the slope angle $\alpha$ is not too large, and we approximate $\cos \alpha(s(t)) \simeq 1$ in condition~(\ref{eqn_no_slip}).
Setting $F(t)= T_f(t)+ T_r(t)$ as the overall traction and braking force, and $F_r(t)=R_f(t)+R_r(t)$ as the overall rolling resistance, we rewrite the first of~\eqref{eqn_dyn} as
\begin{equation}
  \label{eqn_dyn}
M \dot v(t)= -F_a(t) + F(t) - M g \sin \alpha(s(t)) - F_r.
\end{equation}

The aerodynamic drag is $F_a(t)=\frac{1}{2} \rho A_f c_d v(t)^2= \Gamma  v(t)^2$, where $\rho$ is the air density, $A_f$ is the car cross-sectional area, and $c_d$ is the dimensionless drag coefficient.
We model the overall tyres rolling resistance with the simplified model
$F_r(t)= (W_f(t)+W_r(t)) c= M g \cos \alpha(s) c \simeq M g c$, where $c$ is the dimensionless rolling resistance coefficient. Note that $F_r(t)=F_r$ is independent of $t$. The power consumed (if positive) or recovered (if negative) by the vehicle is given by 
\[
P(t)= \max\{\eta F(t),F(t)\} v(t).
\]
According to the value of $\eta \in [0, 1]$, we obtain different power consumption behaviors: in case of a thermal engine, the vehicle is not recovering energy when braking ($\eta = 0$). Otherwise, in case of hybrid or electric engines ($\eta > 0$), the regenerative braking system of the vehicle enables energy recovery during deceleration, which accounts for a negative value for the consumed power, meaning that power is recovered.
\newline\newline\noindent
If we consider as variable the kinetic energy as a function of position $s$, that is $w(s)=\frac{1}{2} M v^2(s)$, essentially it is like considering the squared velocity (see, for instance,~\cite{verscheure09}) scaled by a constant factor. Then, in this way, $w'(s)=M v(s) \dot v(s) v(s)^{-1}=M  \dot v(s)$, and~\eqref{eqn_dyn} becomes
\[
w'(s)= -F_a(s) + F(s) - M g \sin \alpha(s) - F_r,
\]
where the forces are now considered as functions of distance $s$.
We can then write our optimization problem as follows
\begin{align}
& \min_w \int_{s=0}^{s_f} \left(\lambda \max\{\eta F(s),F(s)\} + \frac{1}{\sqrt{w(s)}}\right)ds		& \nonumber \\
\textrm{such that} &											& \nonumber \\
& w(s) \leq \wmax(s)											& \text{for } s \in [0, s_f], \nonumber \\
& |F(s)| \leq M g \mu											& \text{for } s \in [0, s_f], \nonumber \\
& F(s) \leq \frac{\Pmax}{\sqrt{w(s)}}								& \text{for } s \in [0, s_f], \nonumber \\
& w'(s)= -F_a(s) + F(s) - M g \sin \alpha(s)	 - F_r					& \text{for } s \in [0, s_f], \nonumber \\
& w(s) \geq 0												& \text{for } s \in [0, s_f], \nonumber \\
& w(0) = \winit.												& \nonumber
\end{align}
Note that:
\begin{itemize}
\item an upper bound $\wmax(s)$ is imposed for the speed $w(s)$, depending on the position along the path, since, e.g., the maximum allowed speed is different along a curve and along a straight road;
\item an upper bound $\mu M g$ is imposed for the traction or braking force;
\item an upper bound $\Pmax$ is imposed for the power.
\end{itemize}
We discretize the problem by observing variable and input functions at positions $\{0,h,2h,\ldots,(n-1)h\}$,  where $h$ is the discretization step and $n \in \mathbb{N}$ is such that $(n-1)h \approx s_f$ is the length of the path.
We also approximate derivatives with finite differences and the integral in the objective function
with the Riemann sum of the intervals.
Then, the resulting discretized version of the previous problem is:
\begin{align}
& \min_w \sum_{i=1}^{n-1} \left(\lambda \max\{\eta F_i,F_i\} + \frac{1}{\sqrt{w_i}}\right)h		& \nonumber \\
\textrm{such that} &										& \nonumber \\
& w_i \leq \wmax_i										& \text{for } i \in \{1, \ldots, n\}, \nonumber \\
& |F_i| \leq M g \mu										& \text{for } i \in \{1, \ldots, n-1\}, \nonumber \\
& F_i \leq \frac{\Pmax}{\sqrt{w_i}}							& \text{for } i \in \{1, \ldots, n-1\}, \nonumber \\
& \frac{1}{h} (w_{i+1}-w_i) = -\Gamma w_i+ F_i - M g(\sin \alpha_i - c)	& \text{for } i \in \{1, \ldots, n-1\}, \nonumber \\
& w_i \geq 0											& \text{for } i \in \{1, \ldots, n\}, \nonumber \\
& w_1 = \winit,											& \nonumber
\end{align}
which is equivalent to Problem~\eqref{obj_fcn}--\eqref{constr:init_cond}.
%\bibliographystyle{plain} 
%\bibliography{biblio,VelPlan,VelPlan2}

\end{document}